# A highly efficient finite volume method with a diffusion control parameter for hyperbolic problems

Wassim Aboussi[1,2*], Moussa Ziggaf[1,3,4], Imad Kissami[4] and Mohamed Boubekeur[1]

[1*]Sorbonne Paris Nord University, LAGA, CNRS, UMR 7539, Villetaneuse, France.
[2]LAMA, Sidi Mohamed Ben Abdellah University, Faculty of Sciences Dhar El Mahraz, Fez, Morocco.
[3]ENSAO, LMCS, Complexe Universitaire, B.P. 669, 60000 Oujda, Morocco.
[4]MSDA, Mohammed VI Polytechnic University, Lot 660, 43150 Ben Guerir, Morocco.

*Corresponding author(s). E-mail(s): aboussi@math.univ-paris13.fr;
Contributing authors: ziggaf@math.univ-paris13.fr; Imad.kissami@um6p.ma; boubekeur@math.univ-paris13.fr;

**Abstract**

This article proposes a highly accurate, fast and conservative method for hyperbolic systems using the finite volume approach. This innovative scheme constructs the intermediate states at the interfaces of the control volumes using the method of characteristics. The approach is simple to implement, generates entropic solutions, and avoids solving Riemann problems. A diffusion control parameter is introduced to increase the accuracy of the scheme. Numerical examples are presented for the Euler equation for an ideal gas. The results demonstrate the method's ability to capture contact discontinuity and shock wave profiles with high accuracy and low cost as well as its robustness.

**MSC Classification(2020):** 65M08 , 35L65 , 76M12 , 76N15.

# 1 Introduction

A complementary approach for experiment and modeling, numerical simulation is one of the three pillars of scientific research. Fluid mechanics is one of the pioneering sectors in this triptych, and development of numerical schemes well suited to fluid mechanics is a subject that interest numerical scientists. One of the difficulties is reconciling accuracy and robustness with a reasonable computational cost, but the complications can be quite different depending on the targeted applications. Thus, despite the numerous works and the advances in a subject that is still relevant today [1, 2], it is quite natural that there is no uniformly efficient method in all regimes.

In the context of the numerical approximation of hyperbolic systems of conservation laws, several methods based essentially on the solution of the Riemann problem have been retained, and concern shock capturing schemes [3–7]. These methods propose strategies allowing the exact solution of the Riemann problem on each interface, which makes them expensive. In order to reduce the computational time, other approaches propose an approximate solution to the Riemann problem. For example, Roe and Hartan [6, 8, 9] provided schemes based on the evaluation of numerical flux from the exact solution of the linearized problem, and the industry widely uses these schemes because they can capture shock waves with reasonable accuracy.

The main goal of this paper is to describe a new method that can be an excellent tool for simulating most compressible flow phenomena. The proposed scheme belongs to a class of numerical schemes that incorporates the method of characteristics in reconstructing numerical flux. This approach is called the Finite Volume Characteristics (FVC) scheme and has been proposed by Benkhaldoun and Seaïd in [10] and used in the context of shallow water flow [11–14].

The proposed scheme is easy to implement, fast and it accurately solves hyperbolic systems of conservation laws, moreover, it avoids the resolution of Riemann problem in the time integration process, and it is conservative. To approximate the characteristic curves, an iterative process is used, and the intermediate states are calculated using polynomial interpolation. These features are demonstrated using several reference problems for the Euler equations [1, 15]. The presented results provide accurate solutions with a low computational cost. The implementation procedure is described. It is simple to program and generates the numerical results for compressible Euler equations.

This paper is structured as follows: a brief description of the mathematical model will be presented in Section 2. Section 3 is devoted to presenting the process of construction of the numerical scheme. Then, the results of the simulation will be presented in Section 4. The accuracy and the efficiency of the method are discussed and some conclusions are presented in Section 5.

## 2 Governing equation

We consider the one-dimensional Euler equations modeling the dynamics of non-viscous fluid [16, p. 1-12]

$$\frac{\partial \mathbf{W}}{\partial t} + \frac{\partial \mathbf{F}(\mathbf{W})}{\partial x} = 0, \tag{1}$$

where

$$\mathbf{W} = \begin{pmatrix} \rho \\ \rho u \\ E \end{pmatrix}, \quad \text{and} \quad \mathbf{F}(\mathbf{W}) = \begin{pmatrix} \rho u \\ \rho u^2 + p \\ u(E + p) \end{pmatrix}, \tag{2}$$

where $\rho$ is the density of the fluid, $u$ is the fluid particle velocity, $E$ is the total energy and $p$ is the pressure. The state equation links $\rho$, $u$, $p$ and $E$ as

$$E = \rho(\frac{1}{2}u^2 + e(\rho, p)), \tag{3}$$

where $e$ is the specific internal energy; for ideal gases it has a following expression

$$e(\rho, p) = \frac{p}{(\gamma - 1)\rho}, \tag{4}$$

with $\gamma$ is the ratio of specific heats, it is a constant that depends on the particular gas, e.g $\gamma = 1.4$ for air. Another quantity that expresses the ratio of the local velocity of the fluid to the sound speed in this same fluid is called Mach number, which is a dimensionless number defined as

$$M = \frac{u}{c}, \tag{5}$$

where $c$ is the sound speed in a gas. For ideal gas, we can express $c$ by:

$$c = \sqrt{\frac{\gamma p}{\rho}}. \tag{6}$$

$c$ varies with the nature and temperature of the fluid. So, the Mach number does not correspond to a fixed speed, it depends on local conditions of. Generally, the speed is categories by its corresponding regimes [17]. For example, subsonic regime ($M < 0.8$), transsonic regime ($0.8 \leqslant M < 1.2$), supersonic regime ($1.2 \leqslant M < 5$) and hypersonic regime ($M \geqslant 5$).

A simple calculation shows that the system can be written in a quasi-linear form

$$\partial_t \mathbf{W} + \mathbf{J}_F \partial_x \mathbf{W} = 0, \tag{7}$$

where $\mathbf{J}_F$ is the Jacobian matrix of the physical flux $\mathbf{F}$

$$\mathbf{J}_F = \begin{pmatrix} 0 & 1 & 0 \\ (\gamma - 3)\frac{u^2}{2} & (\gamma - 3)u & \gamma - 1 \\ (\frac{\gamma-1}{2}u^2 - H)u & H + (1 - \gamma)u^2 & \gamma u \end{pmatrix},$$

here $H$ is the total specific enthalpy defined by

$$H = \frac{E + p}{\rho}. \tag{8}$$

The matrix $\mathbf{J}_F$ has three eigenvalues: $\lambda_1 = u - c$, $\lambda_2 = u$, and $\lambda_3 = u + c$. and the corresponding right eigenvectors are

$$\mathbf{r}_1 = \begin{pmatrix} 1 \\ u - c \\ H - uc \end{pmatrix}, \quad \mathbf{r}_2 = \begin{pmatrix} 1 \\ u \\ \frac{u^2}{2} \end{pmatrix}, \quad \text{and} \quad \mathbf{r}_3 = \begin{pmatrix} 1 \\ u + c \\ H + uc \end{pmatrix}. \tag{9}$$

There are also other mathematical quantities related to the quasi-linear Euler equations called the Riemann invariants that are constant along the characteristic curves; the right and left Riemann invariants are

$$RI_1 = u - \frac{2c}{\gamma - 1}, \quad \text{and} \quad RI_2 = u + \frac{2c}{\gamma - 1}. \tag{10}$$

Another form of Euler equations using the primitive variables can be formulated as

$$\partial_t \mathbf{W} + \mathbf{A}_{\mathbf{W}} \, \partial_x \mathbf{W} = 0, \tag{11}$$

where

$$\mathbf{W} = \begin{pmatrix} \rho \\ u \\ p \end{pmatrix}, \quad \text{and} \quad \mathbf{A}_{\mathbf{W}} = \begin{pmatrix} u & \rho & 0 \\ 0 & u & \frac{1}{\rho} \\ 0 & \gamma p & u \end{pmatrix}. \tag{12}$$

However, it turns out that for non-smooth solutions the non-conservative formulations give incorrect shock solutions. This point has been noticed for shallow water equations [16, Subsection 3.3]. Despite this, non-conservative formulations have some advantages over their conservative counterpart, when analyzing equations [18, 19].

## 3 Numerical method

In this section, we present our method for Euler equations (1). The method consists of two steps, predictor and corrector; in the first step, we construct the intermediate states using the method of characteristics, while in the second step, the numerical flux in the conservative conservative discretization will be built using the physical flux.

### 3.1 Conservative discretization

The space-time evolution of the fluid motion is described by a vector function $\mathbf{W}(x,t)$, whose components are three flow-dependent variables. Computationally, this function is replaced by $\mathbf{W}_i^n$, which is an approximation of $\mathbf{W}(i\Delta x, n\Delta t)$, where, for simplicity, $\Delta x, \Delta t$ are considered small constants that define a computational space-time grid. To properly capture the shocks generated by the equation system (1), $\mathbf{W}$ was chosen as $\mathbf{W}(x,t)$, the vector of variables in the model continues. P.D. Lax [20] addressed the fundamental problem of determining the $(n+1)$ time-level solution from $n$-level data by creating **interface states** at location such as $(i+1/2)\Delta x$ and proposed the following discretization

$$\frac{\mathbf{W}_i^{n+1} - \mathbf{W}_i^n}{\Delta t} + \frac{\mathbf{F}_{i+1/2}^n - \mathbf{F}_{i-1/2}^n}{\Delta x} = 0. \tag{13}$$

Lax showed that we will recover the partial differential equation (1) when $\Delta x$ and $\Delta t$ go to zero in (13), with a good construction of $\mathbf{F}_{i+1/2}$ from the values of $\mathbf{W}$ in some neighborhood of $i+1/2$. Then, given the initial data $\mathbf{W}_i^0$ for all $i$, we can use equation (13) to construct the solution at the next instant. Thus, we define $\mathbf{W}_i^n$ as follows

$$\mathbf{W}_i^n = \frac{1}{\Delta x} \int_{(i-1/2)\Delta x}^{(i+1/2)\Delta x} \int_{n\Delta t}^{(n+1)\Delta t} \mathbf{W}(x,t) dt dx, \tag{14}$$

and $\mathbf{F}_{i\pm 1/2}^n := \mathbf{F}\left(\mathbf{W}_{i\pm 1/2}^n\right)$ are the numerical fluxes at $x_{i\pm 1/2} := (i \pm 1/2)\Delta x$ and time $t_n := n\Delta t$. The spatial discretization of equation (13) is complete when a numerical construction of the fluxes $\mathbf{F}_{i\pm 1/2}^n$ is chosen. In general, this construction requires a solution of Riemann problems at the interfaces $x_{i\pm 1/2}$. From a computational viewpoint, this procedure is very demanding and may restrict the application of the method for which Riemann solutions are not available.

In the present work, we reconstruct the intermediate states $\mathbf{W}_{i\pm 1/2}^n$ using the method of characteristics. The fundamental idea of this method is to impose a regular grid at the new time level and to backtrack the flow trajectories to the previous time level. At the old time level, the quantities that are needed are evaluated by interpolation from their known values on a regular grid.

### 3.2 Method of characteristics

This method for hyperbolic systems of conservation laws can be carried out componentwise provided that the conservative equations can be rewritten in an advective formulation. In general, the advective form of the system under study is built such that the conservative variables are transported with the same velocity field. In the current study, we apply our method to the Euler equations; we first reformulated the system of equations (1) in an advective form as

$$\partial_t \mathbf{U} + u\partial_x \mathbf{U} = \mathbf{G}(\mathbf{W}), \tag{15}$$

where

$$\mathbf{U} = \begin{pmatrix} \rho \\ \rho u \\ E \end{pmatrix}, \quad \text{and} \quad \mathbf{G}(\mathbf{W}) = \begin{pmatrix} -\rho \partial_x u \\ -\rho \partial_x u - \partial_x p \\ -E \partial_x u - \partial_x (pu) \end{pmatrix}. \tag{16}$$

**Remark 1** *Another formulation using the primitive variables, $\rho$, $u$ and $p$ (see, Eq. (12)) is also possible with the function $\boldsymbol{G}(\boldsymbol{W}) := \left(-\rho\partial_x u, \ -\frac{1}{\rho}\partial_x p, \ \gamma p \partial_x u\right)^T$. We found that for regular solutions this formulation provides the same results as the conservative formulation. However, it is not the same behavior for other solutions, because as already mentioned, this form generates incorrect shocks.*

We calculate now the characteristic curves $x_c(s)$ associated to (15) as

$$\begin{cases} \dfrac{\mathrm{d}x_c(s)}{\mathrm{d}s} = u\left(x_c(s), s\right), \quad s \in \left[t_n, t_n + \alpha_{i+1/2}^n \Delta t\right], \\ x_c\left(t_n + \alpha_{i+1/2}^n \Delta t\right) = x_{i+1/2}, \end{cases} \tag{17}$$

where $u$ is the velocity of the fluid flow. Note that $x_c(s)$ is the departure point at time $s$ of a particle that will arrive at the gridpoint $x_{i+1/2}$ in time $t_n + \alpha_{i+1/2}^n \Delta t$, with $\alpha_{i+1/2}^n$ is a parameter less than 1, that controls the temporal grid, see Fig. 1; the choice of $\alpha_{i+1/2}^n$ is discussed in subsection 3.3. The method of characteristics does not follow the flow particles forward in time, as the Lagrangian schemes do, instead it traces backwards the position at time $t_n$ of particles that will reach the points of a fixed mesh at time $t_n + \alpha_{i+1/2}^n \Delta t$. Therefore, the method avoids the grid distortion difficulties that the conventional Lagrangian schemes have.

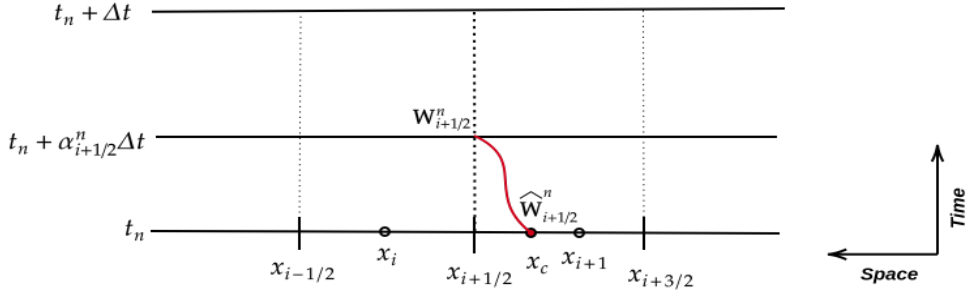

**Fig. 1**: Sketch of the method of characteristics: A fluid particle at gridpoint $x_{i+1/2}$ is traced back in time to $x_c$.

Hence, the solution of (17) can be expressed in an integral form as

$$x_c(t_n) = x_{i+1/2} - \int_{t_n}^{t_n + \alpha_{i+1/2}^n \Delta t} u(x_c(s), s) \, ds. \tag{18}$$

In our simulations we used a first-order Euler method to approximate the integral in (18), however other high-order approximation methods are also possible. In general $x_c(t_n)$ will not coincide with the spatial position of a gridpoint. Thus, once the characteristic curves $x_c(t_n)$ are accurately calculated, the intermediate solutions $\mathbf{W}_{i+1/2}^n$ of a generic function $\mathbf{W}(x_{i+1/2}, t_n)$ are reconstructed using:

$$\mathbf{W}_{i+1/2}^n = \widehat{\mathbf{W}}_{i+1/2}^n + \int_{t_n}^{t_n + \alpha_{i+1/2}^n \Delta t} \mathbf{G}(\mathbf{W}(x_c(s), s)) \, ds, \tag{19}$$

where $\widehat{\mathbf{W}}_{i+1/2}^n = \mathbf{W}(x_c(t_n), t_n)$ are the solutions at the characteristic feet computed by interpolation from the gridpoints of the control volume where the departure points reside, see Fig. 1 for an illustration. For instance, a linear-based interpolation polynomials can be formulated component by component as

$$\widehat{\mathbf{W}}_{i+1/2}^n = \mathbf{W}_i^n + \frac{\mathbf{W}_{i-1}^n - \mathbf{W}_i^n}{\Delta x}(x_c(t_n) - x_i(t_n)), \tag{20}$$

Note that another polynomial interpolation can be used for smooth solutions. However, we have noticed that there is no significant improvement if we change the order of interpolation. This is justified by the fact that the information about the characteristic curve $x_c(t_n)$ which lies between cell $i$ and $i+1$ is given by these cells.

### 3.3 Control parameter $\alpha_{i+1/2}^n$

The choice of the control parameter is based on the stability analysis presented by Benkhaldoun and Seaïd in [10]. This analysis leads us to propose a control parameter $\alpha_{i+1/2}^n$ calculated locally and at each time step with the following formula

$$\alpha_{i+1/2}^n = \tilde{\alpha}_{i+1/2} + \left(\frac{1}{2} - \tilde{\alpha}_{i+1/2}\right) \phi(r_{i+1/2}) \tag{21}$$

where

$$\tilde{\alpha}_{i+1/2} = \frac{\Delta x}{2\Delta t S_{i+1/2}}, \quad \text{and} \quad S_{i+1/2} = \max_k \left(\max_i \left(|\lambda_i^k|, |\lambda_{i+1}^k|\right)\right) \tag{22}$$

here $\lambda_i^k$ is the $k^{th}$ eigenvalue of (7), $S_{i+1/2}$ is the local Rusanov speed and $\phi(r_{i+1/2})$ is a slope limiter. The results presented in Section 4, were obtained using the Minmod limiter; note that other slope limiters functions can be used as van Albada function. The ratio $r_{i+1/2}$ is given by

$$r_{i+1/2} = \frac{q_i - q_{i-1}}{q_{i+1} - q_i}, \tag{23}$$

where

$$q_i = \max_i \left( \left| u_i + \frac{2c_i}{\gamma - 1} \right|, \left| u_i - \frac{2c_i}{\gamma - 1} \right| \right), \tag{24}$$

We note that if $\phi = 1$, we have $\alpha_{i+1/2}^n = \frac{1}{2}$ and as proved in [10, Lemma. 3.2] the method is a second order scheme. On the other hand, when $\phi = 0$, then $\alpha_{i+1/2}^n = \tilde{\alpha}_{i+1/2}$ which, combined with (25), give us a stable and TVD scheme (thanks to Lemma 2.2 [10]).

## 4 Numerical results and discussions

In order to evaluate the accuracy, performance, and robustness of our method, we applied it to Euler equations (2) with a series of numerical tests. In this paper, we consider a perfect gas with a specific heat ratio $\gamma = 1.4$. The computational domain is $[0, 1]$ and the boundary conditions are transmissive. The exact solutions are performed using the open source code [21].

We compare our results with those given by Rusanov [22], Roe [6] and HLL [23] methods. For time discretization, we use a first order Euler scheme and the time step $\Delta t$ is adjusted according to the following CFL condition

$$\Delta t = Cr \frac{\Delta x}{\max_i (\alpha_{i+1/2}^n \Lambda_i^n)} \tag{25}$$

where $Cr$ is the Courant number and $\Lambda_i = \max_k(|\lambda_i^k|)$ is the spectral radius of the Euler equations.

### 4.1 Sod shock tube

Sod shock tube problem [24] with a sonic point in rarefaction is one of the most important tests since it evaluates the satisfaction of the entropy property of numerical methods. The solution of this problem consists of a right shock wave, a right travelling contact discontinuity wave and a left sonic rarefaction wave [16]. The initial condition is defined as

$$\begin{cases} \rho^0(x < 0.5) = 1, \\ u^0(x < 0.5) = 0.75, \\ p^0(x < 0.5) = 1, \end{cases} \text{ and } \begin{cases} \rho^0(x \geq 0.5) = 0.125, \\ u^0(x \geq 0.5) = 0, \\ p^0(x \geq 0.5) = 0.1 \end{cases} \tag{26}$$

The numerical results presented in this section were calculated with 200 cells and Cr= 0.8. In Fig. 2 we compare our method to the exact solution at time $t = 0.2s$.

We note that the shock, the contact and the rarefaction are correctly captured. Fig. 3 shows a comparison between our method, Roe scheme, HLL and Rusanov schemes. It is clear that our method is more accurate. Another essential advantage of the FVC scheme is that it perfectly approximates the rarefaction wave, including the sonic point which is not the case for the Roe scheme where the entropy problem appears. The same error has been observed in [16, p. 227] for the Godunov scheme and other schemes (see [16, p. 280]).

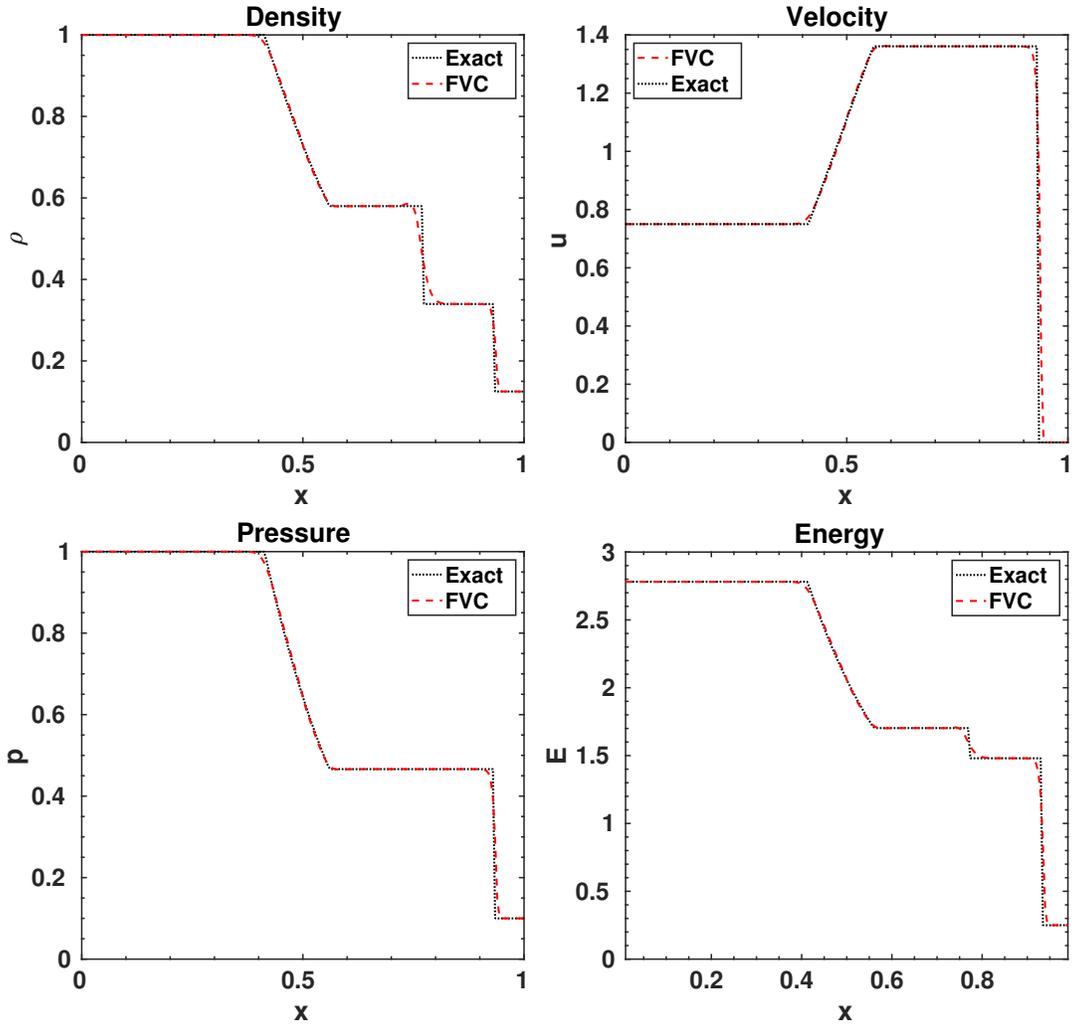

**Fig. 2**: Sod shock tube: density $\rho$ (top left), velocity $u$ (top right), pressure $p$ (bottom left) and total energy $E$ (bottom right) at time $t = 0.2s$ with 200 regular cells.

Fig. 4 represents the behaviour of FVC scheme respect to the choice of $\alpha_{i+1/2}^n$. It is clear that choosing $\alpha_{i+1/2}^n = \frac{1}{2}$ implies the creation of oscillations. On the other hand, choosing $\alpha_{i+1/2}^n = 1$ gives a stable scheme, but the results are more diffused than the case of a variable $\alpha_{i+1/2}^n$. These numerical results are consistent with our comments in section (3) on the reason behind the choice of the control parameter $\alpha_{i+1/2}^n$.

In Fig. 5 we plot the Riemann invariants (left) and the variations of $\alpha_{i+1/2}^n$ computed using Riemann invariants (right), this figure shows us that $\alpha_{i+1/2}^n$ is involved in the area where the shock, the rarefaction and the contact discontinuity appear.

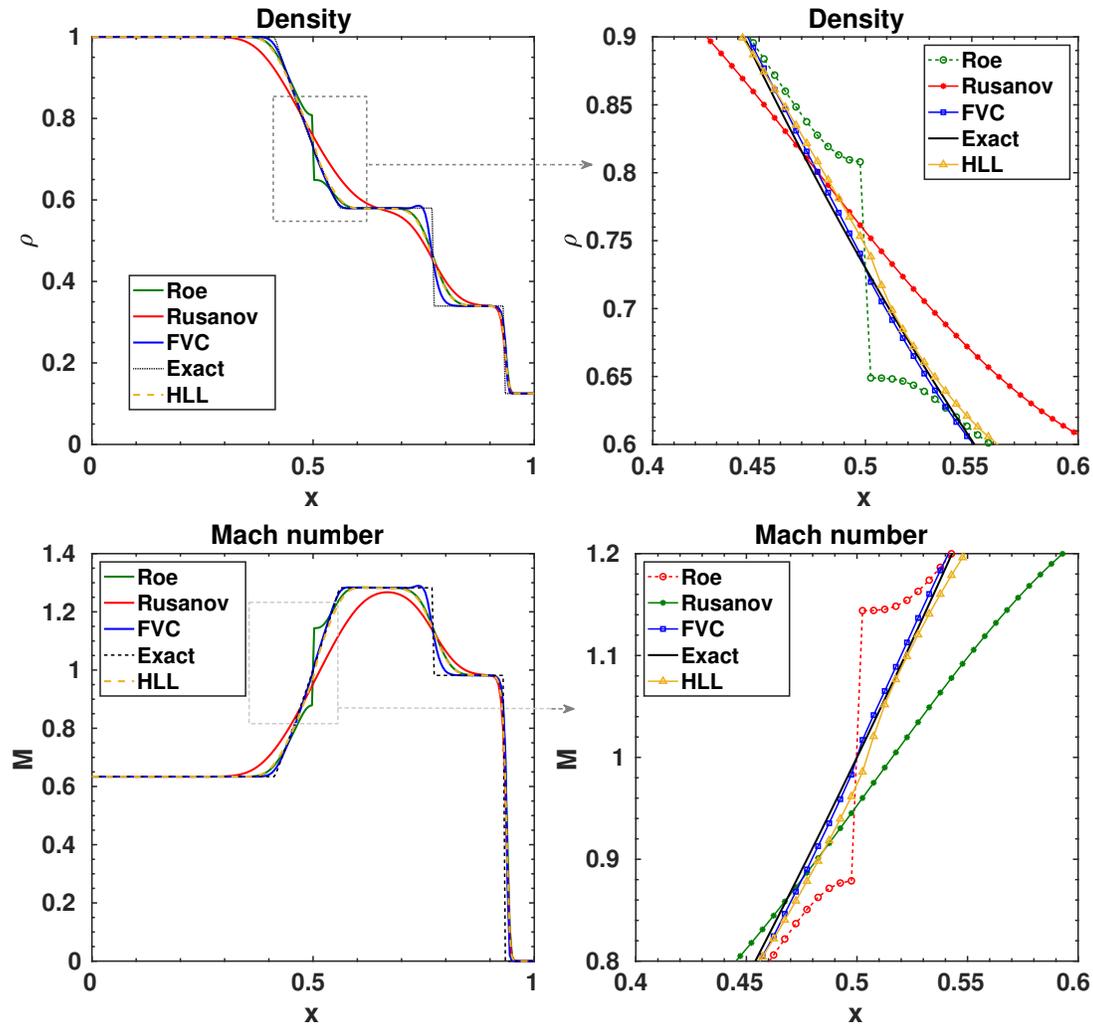

**Fig. 3**: Sod shock tube: density $\rho$ (top left), a zoom on $\rho$ around the sonic point (top right), Much number $M$ (bottom left) and a zoom on $M$ around the sonic point (bottom right) at time $t = 0.2s$ with 200 regular cells.

In Table 1 and Fig. 6, we compare the $L^1$ error for the shock tube problem using different schemes. This table confirms what we said previously about the accuracy of the FVC method compared to the Rusanov scheme, HLL scheme, and a modified Roe scheme. The convergence rate for our method is close to 0.75 while it is equal to 0.61 and 0.56 for HLL and Rusanov respectively.
In table Table 2 we present the computational times for each method. As shown in this table, our method is faster than other schemes, which is one of the most important advantages of our method.

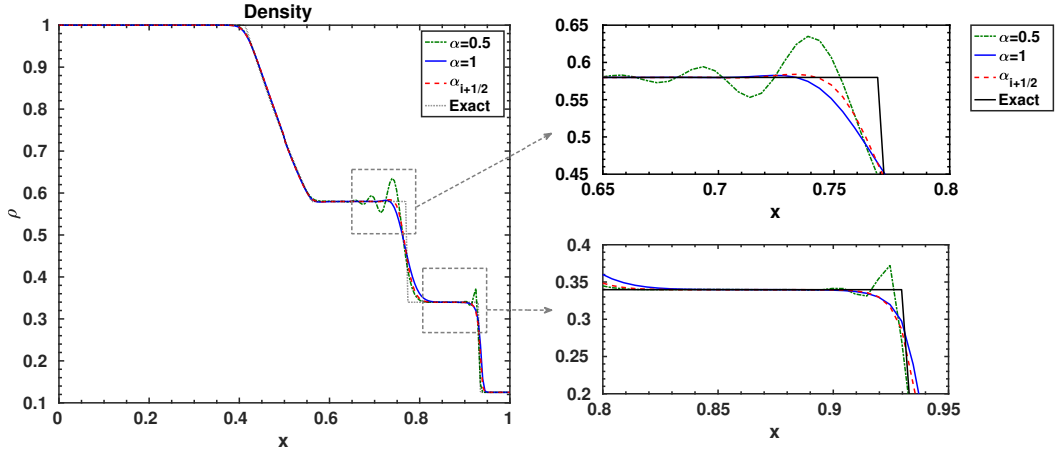

**Fig. 4**: Sod shock tube: numerical solution profile according to the choice of the parameter $\alpha_{i+1/2}^n$ at time $t = 0.2s$ with 200 regular cells.

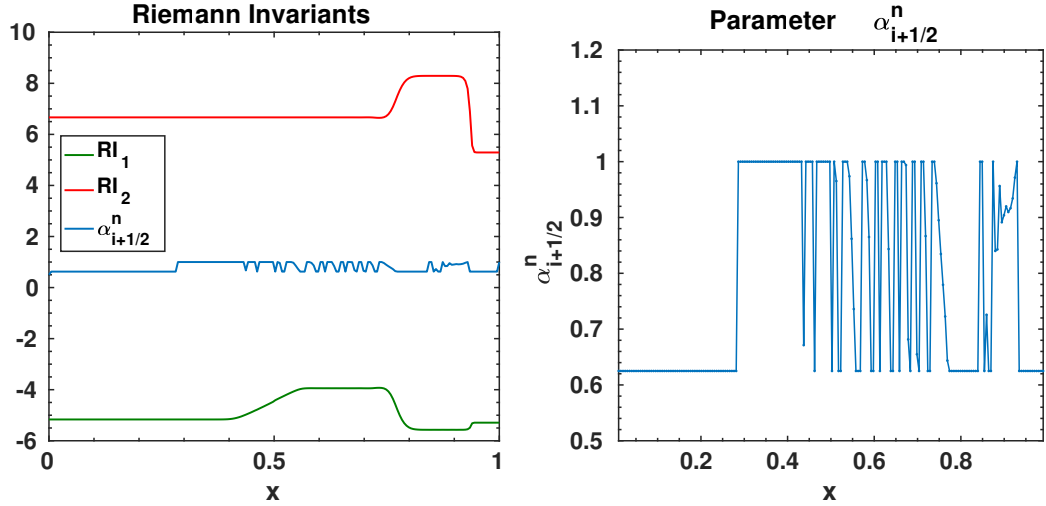

**Fig. 5**: Sod shock tube: Riemann invariant (left) and parameter $\alpha_{i+1/2}^n$ (right) at time $t = 0.2s$ with 200 regular cells.

**Table 1**: Sod shock tube: $L^1$-error for the density at time $t = 0.2s$

| Gridpoints | Rusanov | Roe* | HLL | FVC |
| --- | --- | --- | --- | --- |
| 100 | 3.087960e-02 | 1.549931e-02 | 1.568316e-02 | 7.757252e-03 |
| 200 | 2.148587e-02 | 1.006690e-02 | 1.007645e-02 | 4.421196e-03 |
| 400 | 1.460357e-02 | 6.607027e-03 | 6.665433e-03 | 2.668536e-03 |
| 800 | 9.644334e-03 | 4.372793e-03 | 4.387840e-03 | 1.550064e-03 |
| 1600 | 6.311645e-03 | 2.898692e-03 | 2.900562e-03 | 8.843113e-04 |
| 3200 | 4.121060e-03 | 1.948443e-03 | 1.949273e-03 | 5.363624e-04 |

* Roe with Harten entropy correction [9].

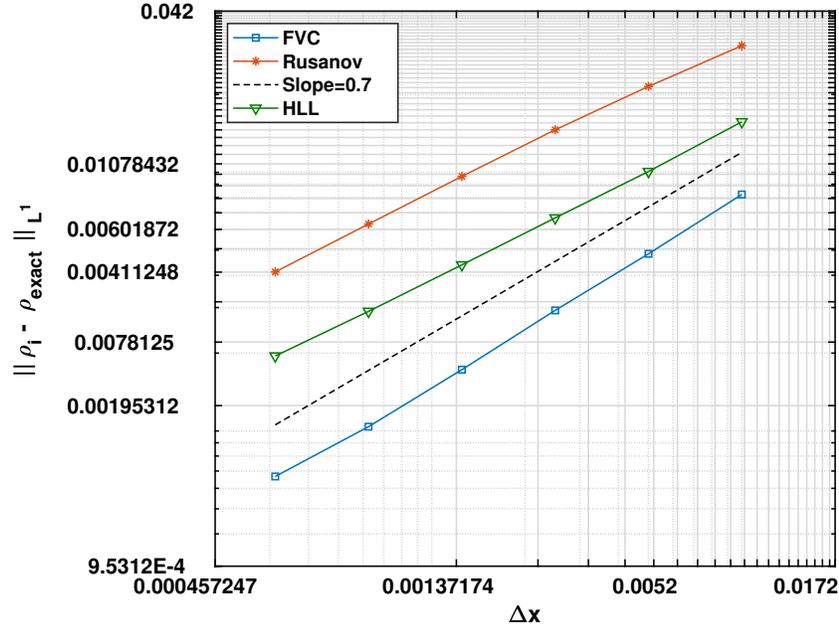

**Fig. 6**: $L^1$ error plot (logarithmic scales) for sod shock problem at time $t = 0.2s$.

**Table 2**: Computational times in seconds for sod shock tube problem.

| Gridpoints | Rusanov | Roe* | HLL | FVC |
| --- | --- | --- | --- | --- |
| 100 | 0.33 | 0.54 | 0.32 | 0.25 |
| 200 | 1.22 | 2.03 | 1.10 | 0.88 |
| 400 | 5.14 | 8.08 | 4.16 | 3.35 |
| 800 | 18.68 | 32.85 | 17.02 | 13.63 |
| 1600 | 80.82 | 142.34 | 71.55 | 57.36 |
| 3200 | 328.79 | 541.69 | 263.22 | 222.68 |

Note: The CPU time was measured on a computer with Intel i7-1165G7 @2.80GHz×8 processor
* Roe with Harten entropy correction [9].

## 4.2 Vacuum test

We now turn to the well-known vacuum test, which is used to evaluate the performance of numerical methods for low density flows. The solution consists of two symmetric rarefaction waves and a stationary contact discontinuity. This problem can be found in [16] and the initial conditions are

$$\begin{cases} \rho^0(x < 0.5) = 1, \\ u^0(x < 0.5) = 0.4, \\ p^0(x < 0.5) = -2.0, \end{cases} \text{ and } \begin{cases} \rho^0(x \geq 0.5) = 1, \\ u^0(x \geq 0.5) = 0.4, \\ p^0(x \geq 0.5) = 2.0. \end{cases} \quad (27)$$

In Fig. 7 we compare the FVC scheme to the numerical solution computed with the Rusanov method, HLL method and to the exact solution with 200 regular cells and Cr=0.8 at $t = 015s$. The rarefaction waves are captured, and we note that in the contact wave zone, where density and pressure are close to zero, the results are acceptable and the positivity of the solution is preserved. We recall that the Roe scheme fails on this problem, but a modified version can be used (see, for example, Einfeldt et al. [25]).

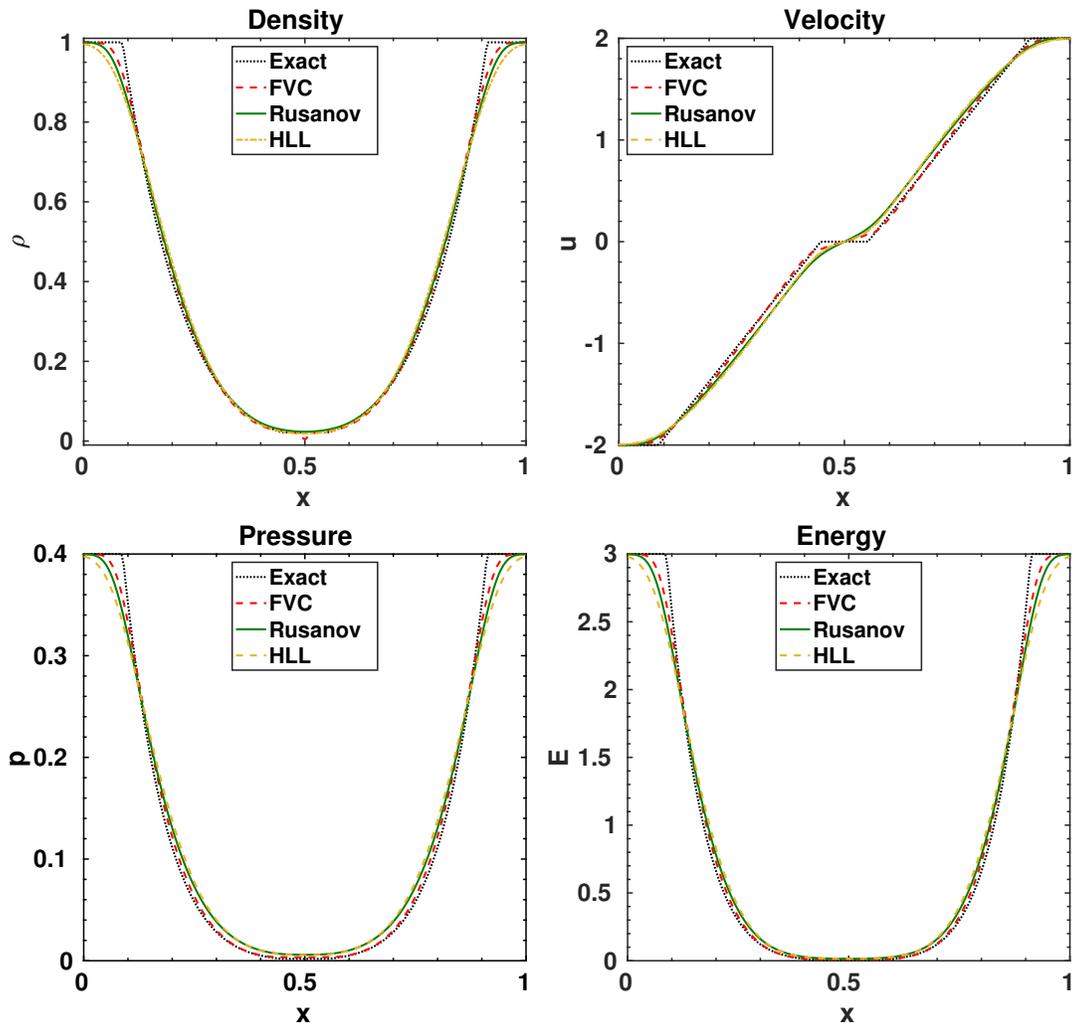

**Fig. 7**: Vacuum test: density $\rho$ (top left), velocity $u$ (top right), pressure $p$ (bottom left) and total energy $E$ (bottom right) at time $t = 0.15s$ with 200 regular cells.

In Fig. 8 we present the variations of $\alpha^n_{i+1/2}$ (right) and the Riemann invariants (left), this figure shows us that $\alpha^n_{i+1/2}$ adapts itself where the stationary discontinuity appears which give us a good approximation of the velocity profile in this zone.

## 4.3 Robustness test

In this section we will perform the robustness of our method using two benchmarks, these tests was proposed in [16]. The first one consists of a strong right shock wave, a contact discontinuity and a left rarefaction wave, and we have:

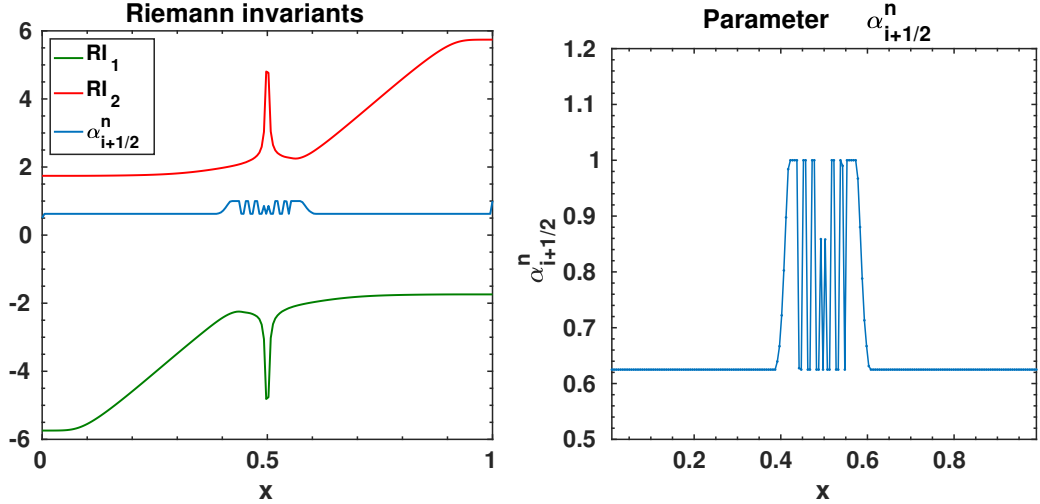

**Fig. 8**: Vacuum test: Riemann invariant (left) and parameter $\alpha^n_{i+1/2}$ (right) at time $t = 0.15$ s with 200 regular cells.

$$\begin{cases} \rho^0(x < 0.5) = 1, \\ u^0(x < 0.5) = 0, \\ p^0(x < 0.5) = 1000, \end{cases} \text{and} \quad \begin{cases} \rho^0(x \geq 0.5) = 1, \\ u^0(x \geq 0.5) = 0, \\ p^0(x \geq 0.5) = 0.01 \cdot \end{cases} \quad (28)$$

the second test consists of a strong left shock wave, a contact discontinuity and a right rarefaction wave; the initial conditions are

$$\begin{cases} \rho^0(x < 0.5) = 1, \\ u^0(x < 0.5) = 0, \\ p^0(x < 0.5) = 0.01, \end{cases} \text{and} \quad \begin{cases} \rho^0(x \geq 0.5) = 1, \\ u^0(x \geq 0.5) = 0, \\ p^0(x \geq 0.5) = 100 \cdot \end{cases} \quad (29)$$

In Fig. 9 we present the numerical solutions for Euler equations with initial data (28) obtained by FVC scheme, Roe scheme, HLL scheme and by Rusanov scheme with 2000 grid cells at time t=0.012s. All schemes show a correct agreement with the exact solution but as we can see in the density curve (left top), FVC scheme is more accurate on the contact discontinuity.

Fig.10 shows the numerical solutions for Euler equations with initial data (29) obtained by FVC scheme, Roe scheme, HLL and by Rusanov scheme with 2000 cells at time $t = 0.035s$.
The difference between this problem and the previous one, is that the velocity is negative. As in Fig. 9, we note that our method captured the contact discontinuity better then other schemes.

### 4.4 A low speed contact discontinuity

In this section we check the ability of our method to resolve slowly–moving contact discontinuities and also a stationary contact discontinuity. Toro et al.

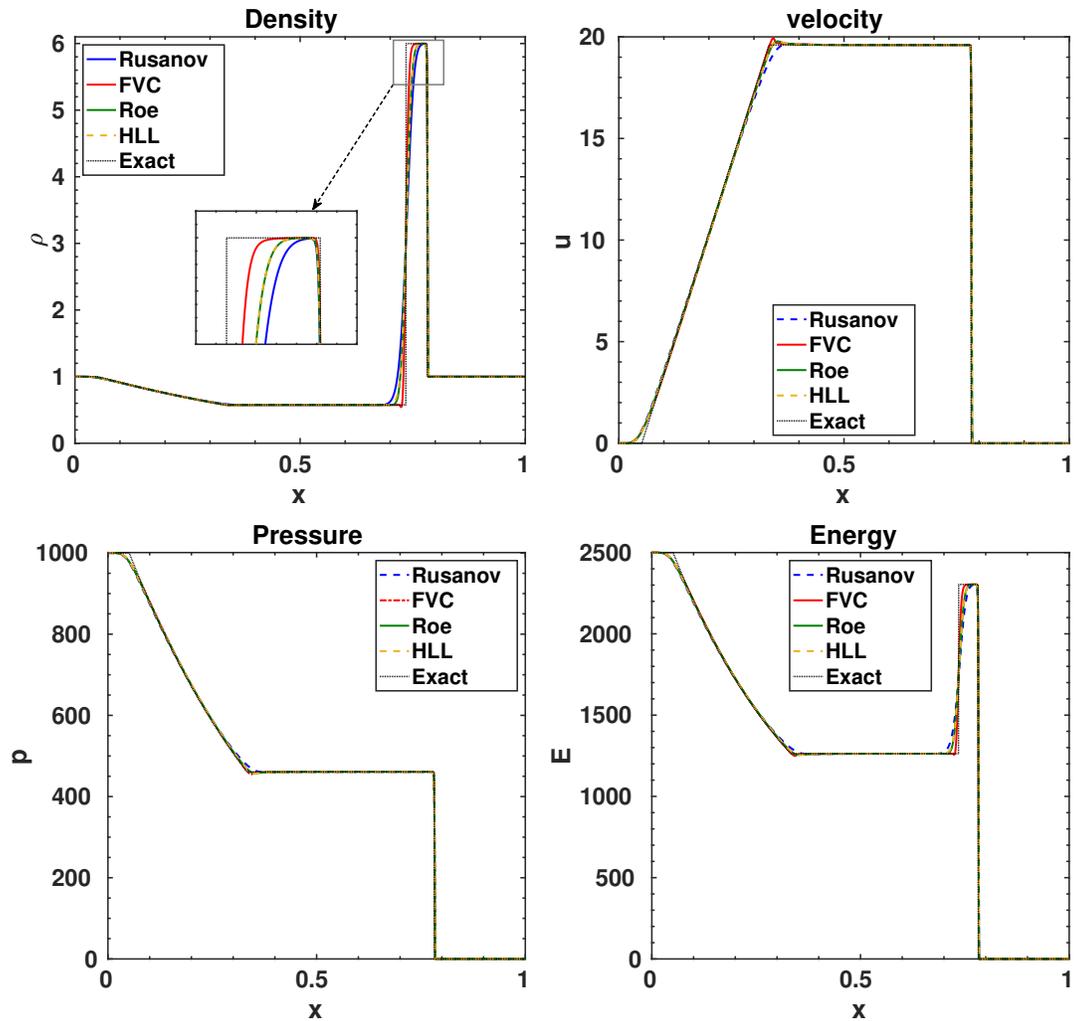

**Fig. 9**: Robustness test: density $\rho$ (top left), velocity $u$ (top right), pressure $p$ (bottom left) and total energy $E$ (bottom right) at time $t = 0.012s$ with 2000 regular cells.

[16] proposed two problems; the first one corresponds to an isolated stationary contact wave and the initial data is given by

$$\begin{cases} \rho^0(x < 0.5) = 1.4, \\ u^0(x < 0.5) = 0, \\ p^0(x < 0.5) = 1, \end{cases} \text{ and } \begin{cases} \rho^0(x \geq 0.5) = 1, \\ u^0(x \geq 0.5) = 0, \\ p^0(x \geq 0.5) = 1. \end{cases} \quad (30)$$

The second one corresponds to an isolated contact moving slowly to the right where the initial data is

$$\begin{cases} \rho^0(x < 0.5) = 1.4, \\ u^0(x < 0.5) = 0.1, \\ p^0(x < 0.5) = 1, \end{cases} \text{ and } \begin{cases} \rho^0(x \geq 0.5) = 1, \\ u^0(x \geq 0.5) = 0.1, \\ p^0(x \geq 0.5) = 1. \end{cases} \quad (31)$$

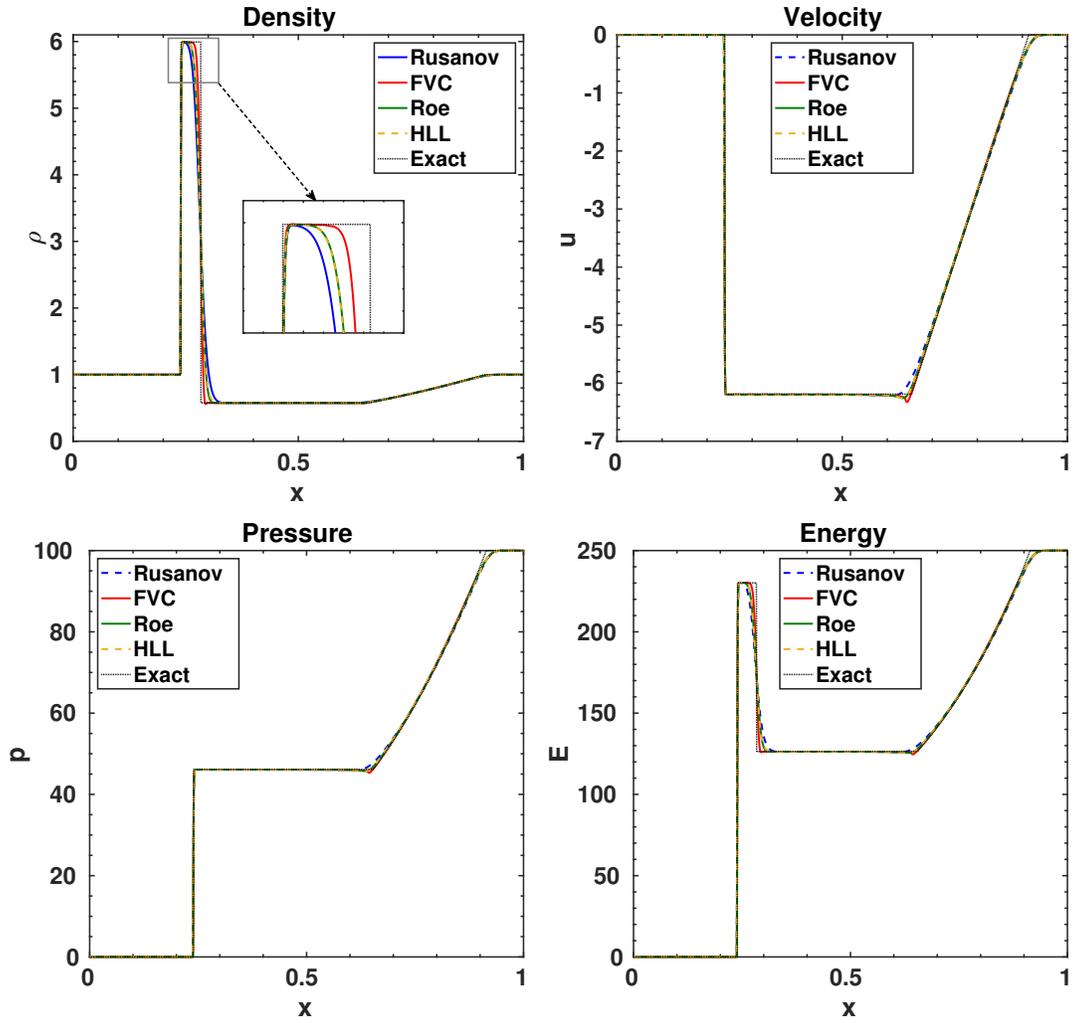

**Fig. 10**: Robustness test: density $\rho$ (top left), velocity $u$ (top right), pressure $p$ (bottom left) and total energy $E$ (bottom right) at time $t = 0.035s$ with 2000 regular cells.

Fig. 11 and Fig. 12 show the numerical results obtained by FVC, HLL, Roe and Rusanov, compared to the exact solution with 200 cells at time t=2s. In Fig. 11 we can see that the numerical results obtained by our method are very similar to those obtained by Roe scheme where the contact is perfectly captured unlike Rusanov scheme and HLL who diffuse. For the slow moving contact test (31), we remark on the density curve that our method is more accurate than Roe, HLL and Rusanov. For the velocity, which is supposed to be constant, there is a small oscillations that appears in all numerical solutions but smaller in the case of FVC we mention that for this sensitive benchmark we used a constant $\alpha_{i+1/2}^n$.

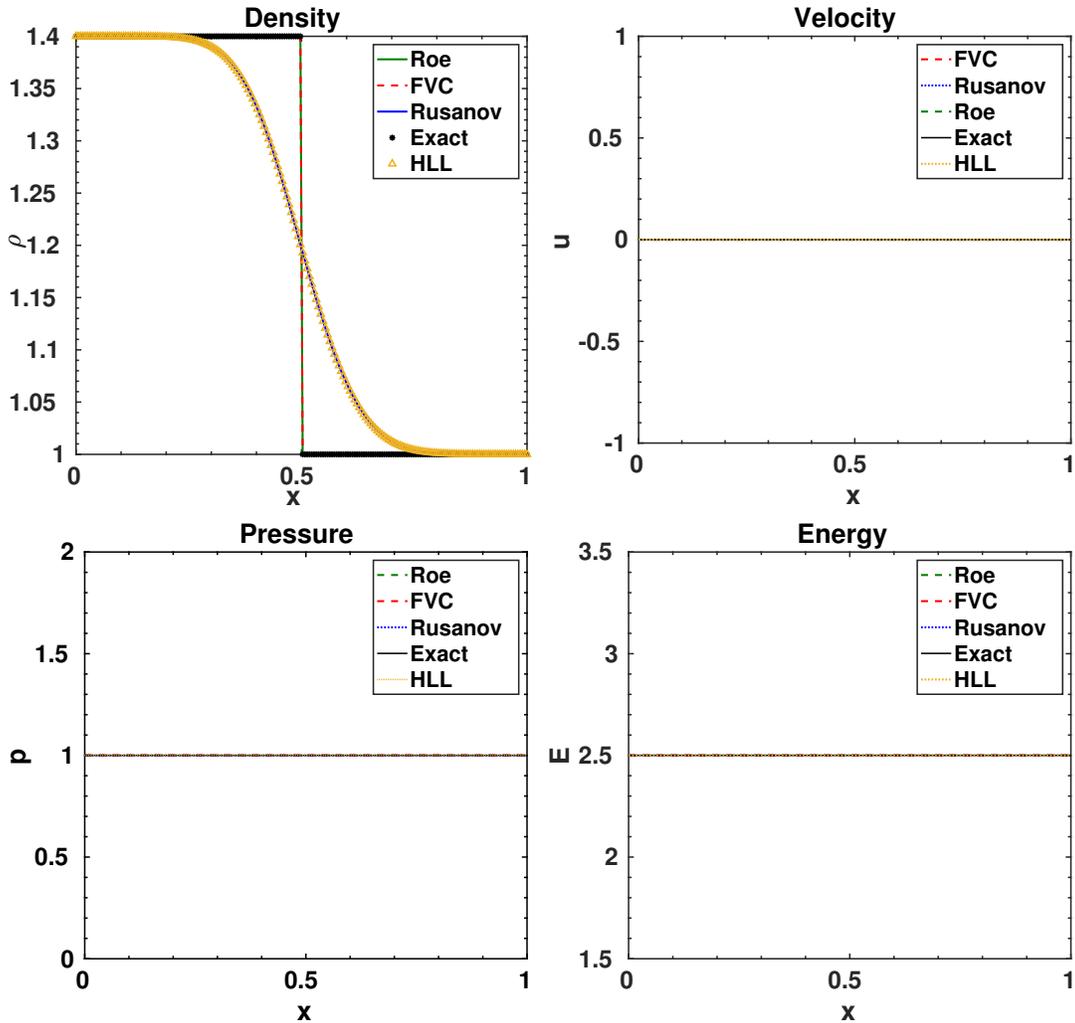

**Fig. 11**: Stationary discontinuity test: density $\rho$ (top left), velocity $u$ (top right), pressure $p$ (bottom left) and total energy $E$ (bottom right) at time $t = 2.0s$ with 200 regular cells.

## 5 Conclusions and outlook

In this work, we proposed an accurate finite volume method for solving hyperbolic problems with application to the one dimensional Euler equations. This method does not need the Jacobian matrix or solving a Riemann problem, which makes it a simple method to implement. The proposed method has been tested using several benchmarks; the results show the high accuracy of our method and more specifically its ability to capture the contact discontinuities. An essential advantage of this method is that it converges to the entropic solution, i.e., the physical solution, without any entropic correction. Moreover the method is fast and highly accurate.

In future works, we will extend this method to multidimensional problems on unstructured meshes, with application to several different physical problems such as flows in pressurized conduits for example.

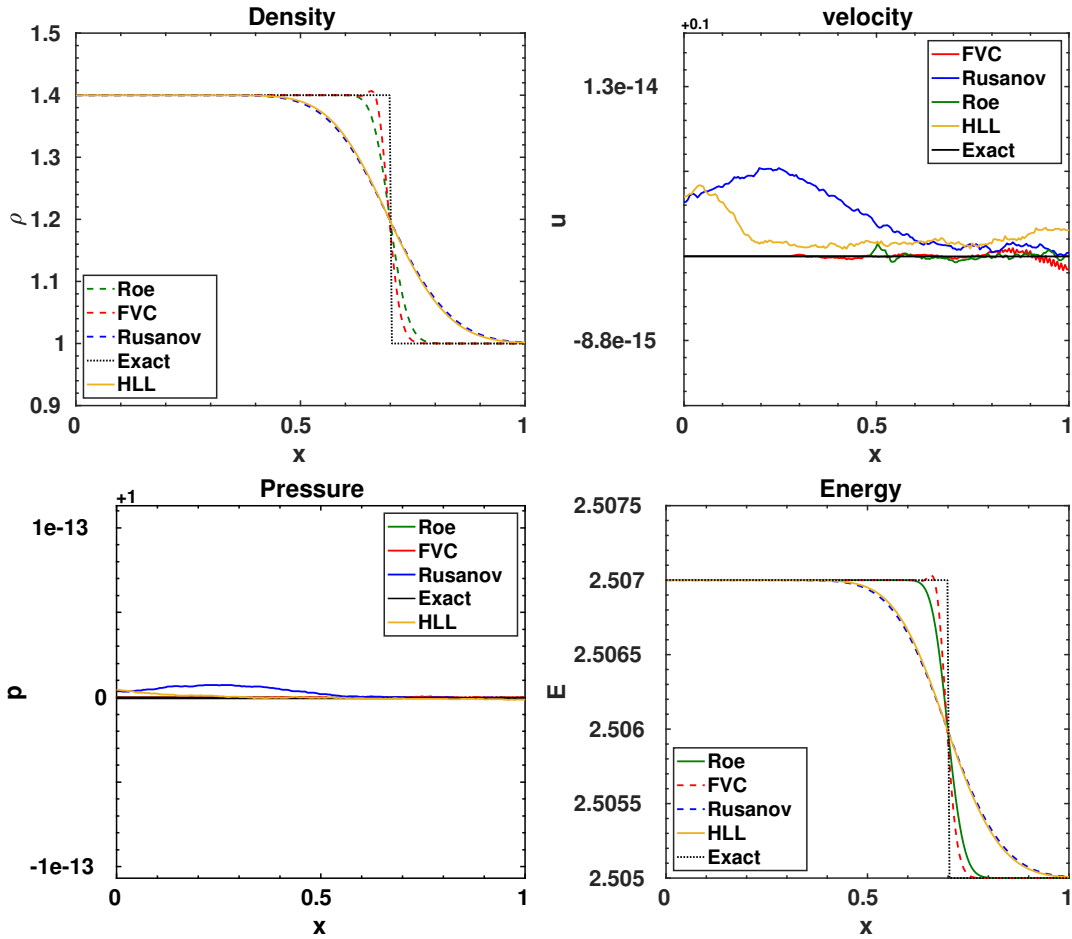

**Fig. 12**: Low speed contact discontinuity test: density $\rho$ (top left), velocity $u$ (top right), pressure $p$ (bottom left) and total energy $E$ (bottom right) at time $t = 2.0s$ with 200 regular cells.

# Acknowledgments

The authors thank Professor F. Benkhaldoun and Professor I. Elmahi for fruitful discussions, orientations and helpful comments.